\documentclass[10pt]{amsart}
\pdfoutput=1
\usepackage{amsmath, amsthm, amssymb,slashed,stmaryrd}
\usepackage{ifpdf}
\usepackage[pdftex]{graphicx}
\usepackage{tikz}
\usetikzlibrary{matrix,arrows,calc}

\usepackage[pdftex,plainpages=false,hypertexnames=false,pdfpagelabels]{hyperref}

\long\def\todo#1{{\color{red} {#1}}}
\long\def\stodo#1{{\color{blue} {#1}}}

 \theoremstyle{plain}
 \newtheorem{thm}{Theorem}[section]

  \newtheorem{quest}{Question}[section]
 \newtheorem{cor}[thm]{Corollary}
 \newtheorem{lem}[thm]{Lemma}
  
 \newtheorem{prop}[thm]{Proposition}
 \theoremstyle{definition}
 \newtheorem{defn}[thm]{Definition}

 \newtheorem{ex}[thm]{Example}
 \theoremstyle{remark}
 \newtheorem{rmk}[thm]{Remark}

\def\beq{\begin{eqnarray}}
\def\eeq{\end{eqnarray}}

\DeclareSymbolFont{bbold}{U}{bbold}{m}{n}
\DeclareSymbolFontAlphabet{\mathbbold}{bbold}
\def\one{\mathbbold{1}}
 \newcommand{\bp}{\begin{proof}[Proof]}
 \newcommand{\ep}{\end{proof}}
\DeclareMathOperator{\SM}{\underline{\sf SMfld}}

\def\pt{\rm pt}

\def\CS{{\rm CS}}

\def\ev{{\rm ev}}
\def\odd{{\rm odd}}

\newcommand{\sq}{\mathord{/\!\!/}}

\def\R{{\mathbb{R}}}
\def\A{{\mathbb{A}}}

\def\id{{{\rm id}}}

\def\proj{{\rm proj}}

\def\K{{\rm {K}}}
\def\C{{\mathbb{C}}}

\def\Z{{\mathbb{Z}}}

\def\Vect{{\sf Vect}} 
\def\Alg{{\sf Alg}} 
\def\pt{{\rm pt}}

\def\Hom{\mathop{\sf Hom}}

\def\SM{ {\underline{\sf SMfld}}}
\def\EFT{ \hbox{-{\sf EFT}}_{\rm eff}}

\def\Path{ {{\sf P}_0}}
\def\Loop{ {{\sf L}_0}}
\def\path{ {{\mathfrak{p}}_0}}

\def\twocommute{\ensuremath{\rotatebox[origin=c]{30}{$\Rightarrow$}}}

\newcommand\triplerightarrow{
\mathrel{\substack{\textstyle\rightarrow\\[-0.5ex]
                      \textstyle\rightarrow \\[-0.5ex]
                      \textstyle\rightarrow}}
}

\begin{document}

\title[Equivariant K-theory and gauged mechanics]{Twisted equivariant differential K-theory from gauged supersymmetric mechanics}

\author{Daniel Berwick-Evans}

\address{Department of Mathematics, University of Illinois at Urbana-Champaign}
\email{danbe@illinois.edu}

\date{\today}
\begin{abstract}
We use the geometry of the space of fields for gauged supersymmetric mechanics to construct the twisted differential equivariant K-theory of a manifold with an action by a finite group. 
\end{abstract}
\maketitle 

\section{Introduction and statement of results}

Let $G$ be a finite group acting on a smooth manifold~$X$ and $\beta\colon G\times G\to U(1)$ a 2-cocycle. Below we construct a model for the $\beta$-twisted equivariant differential K-theory of~$X$ using the geometry of $1|1$-dimensional super paths in the stack~$X\sq G$, i.e., fields for gauged~$N=1$ supersymmetric quantum mechanics. Specifically, we define a flavor of Atiyah--Segal functorial field theory for a bordism category of \emph{energy zero gauged super paths} in~$X$. Ideas from Wilsonian effective field theory motivate a version of stable isomorphism between these functorial field theories, and equivalence classes define differential K-theory. 



\begin{thm}
Let $G$ be a finite group and $\beta\colon G\times G\to U(1)$ a normalized 2-cocycle. There is a natural isomorphism
$$
\widehat{\K}_G^\beta(X)\cong 1|1\EFT^\beta(X\sq G)/{\sim}
$$
between the $\beta$-twisted equivariant differential K-theory of~$X$ and stable isomorphism classes of $\beta$-twisted gauged low-energy effective Euclidean field theories over~$X$.\label{thm1}
\end{thm}


The relationship between 1-dimensional physics and K-theory has been studied from many perspectives. Our approach is heavily influenced by Stolz and Teichner's construction of a \emph{space} of $1|1$-dimensional field theories classifying K-theory~\cite{ST04,HST}, and their subsequent variations on $1|1$-dimensional field theories over~$X$~\cite{ST11} that are expected to construct the K-theory of $X$;~see Remark~\ref{rmk:ST} for a comparison to our model. One point to emphasize is that the cocycles in Theorem~\ref{thm1} paint a different geometric picture when compared to the classifying space description, especially in equivariant refinements. For example, a character map on cocycles leads immediately to the equivariant Chern character as a differential form on fixed point sets for the $G$-action, whereas such a map looks rather unmotivated from the classifying space perspective and requires additional choices for its construction. 


The relationship between \emph{gauged} 1-dimensional physics and \emph{equivariant} K-theory goes back to Alvarez-Gaum\'e's proof of the equivariant index theorem~\cite{Alvarez} (see also Witten~\cite{WitteninStrings}), and our approach is in keeping with these ideas. We have also taken cues from the orbifold K-theory of Adem and Ruan~\cite{AdemRuan}, especially in regard to the equivariant Chern character. 

Another source of inspiration for our methods is Kitaev's~\cite{Kitaev} study of phases of condensed matter systems. He explains how varying the energy cutoff in an effective field theory has the effect of stabilizing vector spaces of states. In the presence of certain (super) symmetries, stable equivalence classes can be identified with K-theory classes. These considerations lead to the definition of stable isomorphism in Theorem~\ref{thm1}.




This paper fits into a larger goal of clarifying relationships among supersymmetric field theories, homotopy theory, and representation theory. A driving force is Witten and Segal's vision \cite{Witten_Elliptic,WittenDirac,Segal_Elliptic} (furthered by Stolz and Teichner~\cite{ST04}) for studying elliptic cohomology and loop group representations using $2|1$-dimensional gauged field theories. One goal below is to provide a framework that permits a leap from $1|1$-dimensional gauged theories to the more complicated~$2|1$-dimensional ones, at least in the case of finite groups. 



\subsection{Notation and terminology}\label{sec:notterm}

Throughout,~$X$ will denote a compact (ordinary) manifold, often regarded as a supermanifold. Our super manifolds have complex algebras of functions, called \emph{cs-manifolds} in Deligne and Morgan's review article~\cite{DM}; otherwise we follow their conventions. We use the functor of points, reserving the letter~$S$ for a test super manifold. A \emph{vector bundle} on a supermanifold is a finitely generated projective module over the ring of functions, and we will use the sections notation, e.g., $\Gamma(E)$, when emphasizing this module structure. A \emph{super connection} on a (super) vector bundle~$E$ over~$X$ will mean a super connection in the sense of Quillen~\cite{Quillensuper}: an odd operator $\A\colon \Omega^\bullet(X;E)\to \Omega^\bullet(X;E)$ satisfying a Leibniz rule over $\Omega^\bullet(X)$. 

For a $G$-action on~$X$, $X\sq G$ denotes the action Lie groupoid or quotient stack; a map $S\to X\sq G$ is a principal $G$-bundle $P$ over~$S$ and a $G$-equivariant map~$P\to X$. Isomorphisms between $G$-bundles give rise to isomorphisms between such maps. 

A \emph{Lie category} $\mathcal{C}$ is a category internal to supermanifolds, meaning we have a supermanifold ${\rm Ob}(\mathcal{C})$ of objects, a supermanifold ${\rm Mor}(\mathcal{C})$ of morphisms, and maps 
$$
{\rm Mor}(\mathcal{C})\times_{{\rm Ob}(\mathcal{C})} {\rm Mor}(\mathcal{C})\triplerightarrow {\rm Mor}(\mathcal{C})\rightrightarrows {\rm Ob}(\mathcal{C})
$$
of supermanifolds encoding source, target, composition, unit and identities (though we have omitted the arrows coming from identity maps above). We typically denote source and target maps by ${\sf s}$ and ${\sf t}$, respectively.

A \emph{representation} of a Lie category is a smooth functor~$\mathcal{C}\to \Vect$: the data is a smooth vector bundle~$E$ over~${\rm Ob}(\mathcal{C})$ and a map of vector bundles~${\sf s}^*E\to {\sf t}^*E$ over~${\rm Mor}(\mathcal{C})$. We require that these data satisfy the obvious compatibility conditions with respect to units and composition. An equivalent definition comes from the functor of points: to an $S$-point of~${\rm Ob}(\mathcal{C})$ we require a vector bundle over~$S$ and to an $S$-point of ${\rm Mor}(\mathcal{C})$ we require a morphism of vector bundles over~$S$. These assignments must be natural in~$S$, and we recover the previous description by taking~$S={\rm Ob}(\mathcal{C})$ and $S={\rm Mor}(\mathcal{C})$.

We also consider \emph{twisted} representations. Let $\Alg^\times$ denote the bicategory of Morita invertible super algebras, invertible bimodules, and bimodule isomorphisms. In brief, a \emph{twist} is a smooth functor~${\sf T}\colon \mathcal{C}\to \Alg^\times$ and a twisted representation is a natural transformation~$\one\Rightarrow{\sf T}$. In more detail, a twist is the data of an algebra bundle~$A$ over ${\rm Ob}(\mathcal{C})$; an invertible ${\sf s}^*A$-${\sf t}^*A$ bimodule bundle~$B$ over ${\rm Mor}(\mathcal{C})$; and a bimodule bundle isomorphism $\beta\colon p_1^*B\otimes p_2^*B\to {\sf c}^*B$ over ${\rm Mor}(\mathcal{C})\times_{{\rm Ob}(\mathcal{C})} {\rm Mor}(\mathcal{C})$, where 
$$
p_1,p_2,{\sf c}\colon {\rm Mor}(\mathcal{C})\times_{{\rm Ob}(\mathcal{C})} {\rm Mor}(\mathcal{C})\to {\rm Mor}(\mathcal{C})
$$
are the projection and composition maps. These satisfy a compatibility encoding associativity over the 3-fold fibered product. In our examples,~$A$ will be the trivial bundle of algebras with fiber~$\C$ and~$B$ will the trivial bundle of bimodules with fiber~$\C$. This leads to twisted representation having the same data as an ordinary representation, but composition in $\mathcal{C}$ is no longer strictly compatible with composition of linear maps: the failure is measured by the isomorphism~$\beta$. 
%
%
\subsection{Outline} In the next section, we define the category of energy zero gauged super paths. In Section~\ref{sec:rep}, we explain how representations of this category define equivariant vector bundles with equivariant super connection, which is the key geometric construction in the paper. Lastly, in Section~\ref{sec:diffK} we prove Theorem~\ref{thm1}.

\section{The category of energy zero gauged super paths}

An \emph{($S$-family of) super paths} in $X$ is a triple $(\gamma,{\sf in},{\sf out})$ for
$$
\gamma\colon S\times \R^{1|1}\to X,\quad {\sf in},{\sf out}\colon S\times \R^{0|1}\hookrightarrow S\times \R^{1|1}\stackrel{\gamma}{\to} X.
$$ 
We call the maps ${\sf in},{\sf out}\colon S\times \R^{0|1}\to X$ the \emph{source} and \emph{target} super points of the super path. There is a super-translation action on super paths from 
$$
\R^{1|1}(S)\times \R^{1|1}(S)\to \R^{1|1}(S),\quad (t,\theta)\cdot(s,\eta)=(t+s+\theta\eta,\theta+\eta), \ (t,\theta),(s,\eta)\in \R^{1|1}(S).
$$
If we can translate the target super point of a given super path to the source super point of another and the paths agree in a neighborhood of this super point, there is an evident concatenation (or composition) of super paths in~$X$. 

When $X$ has a $G$-action, we can consider an analogous setup over $X\sq G$; a map $S\times \R^{1|1}\to X\sq G$ is the data of a principal $G$-bundle $P\to S\times \R^{1|1}$ and a $G$-equivariant map $P\to X$. These maps to $X\sq G$ along with source and target data are fields for the gauged super particle (e.g., see~\cite{WitteninStrings}). We will study a finite-dimensional subspace of these fields that can be given the structure of a Lie category. 

To simplify the discussion, we first give an ad hoc description of the Lie category, and then we explain its relation to the gauged super particle. Define an action of~$\R^{1|1}\times G$ on~$\pi TX$ that factors through the projection homomorphism $\R^{1|1}\times G\to \R^{0|1}\times G$ and then is determined by the action of~$\R^{0|1}$ on~$\pi TX\cong \SM(\R^{0|1},X)$ generated by the de~Rham operator and the $G$-action on~$X$ lifted to differential forms. 

\begin{defn}
Define the Lie category of energy zero gauged superpaths in~$X$, 
$$
\Path(X\sq G):=\left(\begin{array}{c} 
\R_{\ge 0}^{1|1}\times G\times \SM(\R^{0|1},X)\\
{\sf s}\downarrow \downarrow {\sf t}\\
\SM(\R^{0|1},X)
\end{array} \right)\subset \left(\begin{array}{c} 
\R^{1|1}\times G\times \SM(\R^{0|1},X)\\
{\sf s}\downarrow \downarrow {\sf t}\\
\SM(\R^{0|1},X)
\end{array} \right)
$$
as the Lie subcategory of the action Lie groupoid $\pi TX\sq (\R^{1|1}\times G)$. 
\end{defn}

Now we justify the name of this Lie category. An $S$-point of the objects of $\Path(X\sq G)$ is a map $S\times \R^{0|1}\to X$. Using the $G$-action on $X$, we can promote this to the data of a trivial $G$-bundle over $S\times \R^{0|1}$ and a $G$-equivariant map to~$X$, i.e., a map $S\times\R^{0|1}\to X\sq G$. An $S$-point of the morphisms of $\Path(X\sq G)$ is a triple: $g\colon S\to G$, $(t,\theta)\colon S\to \R^{1|1}_{\ge 0}$ and $\phi\colon S\times \R^{0|1}\to X$. From the datum~$\phi$ we build the commutative diagram,
\beq
\begin{array}{c}
\begin{tikzpicture}[node distance=3.5cm,auto]
  \node (A) {$P=S\times G\times \R^{1|1}$};
  \node (B) [node distance= 8cm, right of=A] {$X$};
  \node (C) [node distance = 1.5cm, below of=A] {$S\times \R^{1|1}$};
  \node (D) [node distance = 1.5 cm, below of =B] {$S\times \R^{0|1}$};
  \node (E) [node distance = 4cm, right of =A] {$S\times G\times \R^{0|1}$};
  \draw[->] (A) to node {$\id\times \proj$} (E);
    \draw[->,bend left=20] (A) to node {$\tilde{\phi} $} (B);
    \draw[->] (E) to (B);
  \draw[->] (A) to node {$p$}  (C);
    \draw[->] (D) to [swap] node {$\phi$} (B);
    \draw[->] (C) to [swap] node {$\proj$} (D);
\end{tikzpicture}\end{array}\nonumber
\eeq
where $\proj\colon \R^{1|1}\to \R^{0|1}$ is the projection, and $\tilde{\phi}$ is the unique lift of $\phi\circ \proj$ to a $G$-equivariant map factoring through $\id\times \proj$. By construction, this is a principal $G$-bundle over $S\times \R^{1|1}$ with a $G$-equivariant map to~$X$, i.e., an $S$-family of super paths in~$X\sq G$. We define the source of this super path in~$X\sq G$ as precomposition with the inclusion $S\times G\times \R^{0|1}\subset P\to X$ induced by the standard inclusion~$\R^{0|1}\subset \R^{1|1}$. Define the target map associated to $(t,\theta)\in \R^{1|1}_{\ge 0}(S)$ and $g\in G(S)$ as the composition
$$
S\times G\times \R^{0|1}\subset S\times G\times \R^{1|1}\stackrel{g\times {T}_{(t,\theta)}}{\longrightarrow} S\times G\times \R^{1|1}=P\stackrel{\tilde{\phi}}{\to} X
$$
where the first inclusion is again induced by $\R^{0|1}\subset \R^{1|1}$, and ${ T}_{(t,\theta)}$ is super translation by $(t,\theta)\in \R^{1|1}_{\ge 0}(S)$ and the map~$g\colon S\times G\to S\times G$ comes from left multiplication by~$g\in G(S)$. 


This discussion identifies~$\Path(X\sq G)$ as a subspace of fields of the gauged super particle. These fields also play a preferred role: they have energy zero with respect to the classical action. We first explain the energy zero condition when $G=\{e\}$. Identify a super path $S\times \R^{1|1}\to X$, with an ordinary path $x\colon S\times \R\to X$ and a section $\psi\in \Gamma(x^*\pi TX)$ along the path. This super path has energy zero if~$x$ does in the usual sense of Riemannian geometry, and if $\psi$ is covariantly constant along the path. Hence, we require that~$x$ be a constant path and $\psi$ be a section of $\pi TX$ at the point in~$X$ defined by~$x$. This is the same as a factorization condition on the super path defined by $(x,\psi)$:
$$
S\times \R^{1|1}\to S\times \R^{0|1}\to X. 
$$
As for gauged super paths, because $G$ is finite any $G$-bundle is automatically flat; this means these bundle don't contribute to the energy, so the condition remains essentially unchanged and we require a factorization as above for super paths in $X\sq G$.

It will be useful to observe that energy zero gauged super paths have an evident naturality: if $X\to Y$ is a map of $G$-manifolds, we obtain an induced functor $\Path(X\sq G)\to \Path(Y\sq G)$.

\section{Twisted effective field theories}\label{sec:rep}

One of the primary pieces of data of a twisted effective field theory will be a \emph{twisted representation} of $\Path(X\sq G)$, meaning a smooth natural transformation ${\sf E}$
\beq
\begin{tikzpicture}[node distance=2cm, auto]
  \node at (0,0) (A) {$\Path(X\sq G)$};
  \node at (2.15,0) (B)  {$\Downarrow {\sf E}$ };
  \node at (4,0) (C)  {${\sf Alg}^\times,$};
\draw[->,bend left] (A.20) to node {$\one$} (C.160);
\draw[->, bend right] (A.340) to node [swap] {${\sf T}$} (C.200);
\end{tikzpicture}\label{eq:twist}
\eeq
for ${\sf T}\colon \Path(X\sq G) \to \Alg^\times$ a smooth functor called a \emph{twist} (see Section~\ref{sec:notterm}). 

\begin{rmk} 
Our terminology follows Stolz and Teichner's~\cite{ST11}: the category $\Path(X\sq G)$ is naturally a subcategory of their $1|1$-dimensional super Euclidean bordism category over~$X\sq G$, and (finite-dimensional) twisted field theories in their sense restrict to twisted representations of $\Path(X\sq G)$. 
\end{rmk}

\subsection{Twists from 2-cocycles on $G$}\label{sec:twists}

For any normalized\footnote{A 2-cocycle $\beta$ is \emph{normalized} if $\beta(g,e)=\beta(e,g)=1$ for all $g\in G$. Any 2-cocycle is cohomologous to a normalized one.}  2-cocycle $\beta\colon G\times G\to U(1)$, we define a smooth functor
$$
\Path(X\sq G)\to \Alg
$$
whose bundle of algebras over ${\rm Ob}(\Path(X\sq G))$ is the trivial bundle $A:=\underline{\C}$, and whose bundle of bimodules over ${\rm Mor}(\Path(X\sq G))$ is the trivial bimodule $B:=\underline{\C}$. Over the fibered product 
$$
{\rm Mor}(\Path(X\sq G))\times_{{\rm Ob}(\Path(X\sq G))} {\rm Mor}(\Path(X\sq G))\cong \R^{1|1}_{\ge 0}\times \R^{1|1}_{\ge 0}\times G\times G\times \pi TX,
$$
we take the isomorphism of bimodules determined by~$\beta$,
$$
p^*\beta\colon (p_1^*B)\otimes_A (p_2^*B)\to m^*B,\quad p\colon \R^{1|1}_{\ge 0}\times \R^{1|1}_{\ge 0}\times G\times G\times \pi TX\to G\times G
$$
which is compatible with the units because $\beta$ is normalized and associative because~$\beta$ is a cocycle. 

\begin{rmk} The twist above is the truncation of a fully extended 2-dimensional \emph{topological} field theory called (classical) Dijkgraaf--Witten theory or (classical) Yang--Mills theory. \end{rmk}

\subsection{Twisted representations}
We will characterize $\beta$-twisted representations of $\Path(X\sq G)$ in terms of $\beta$-twisted equivariant vector bundles on~$X$.

\begin{defn} 
A \emph{$\beta$-twisted equivariant vector bundle on $X$} is a vector bundle $E\to X$ and isomorphisms $\rho_x(g)\colon E_x\to E_{gx}$ for each $g\in G$ that vary smoothly with~$x\in X$ so that $\rho_x(e)=\id$ and $\rho_{hx}(g)\circ\rho_x(h)=\beta(g,h)\rho_x(gh)\colon E_x\to E_{ghx}$. An \emph{equivariant super connection} is a super connection $\A$ on $E$ such that $\rho(g)\A =\A\rho(g)$.\end{defn}

\begin{ex} A $\beta$-twisted equivariant vector bundle on the point is a $\beta$-projective representation, i.e., a map $\rho\colon G\to {\rm End}(V)$ such that~$\rho(g)\rho(h)=\beta(g,h)\rho(gh)$ and~$\rho(e)={\rm id}_V$. \end{ex}

\begin{prop} The category of $\beta$-twisted representations of $\Path(X\sq G)$ is equivalent to the category whose objects are $\beta$-twisted equivariant vector bundles on $X$ with equivariant super connection and whose isomorphisms are isomorphisms of bundles over $\pi TX$ compatible with the $G$-action and connection.  \label{prop:geochar}\end{prop}

\bp
The value of a representation on objects is a super vector bundle over $\pi TX$. By trivializing along the fibers of the projection $p\colon \pi TX\to X$, the groupoid of super vector bundles on $\pi TX$ is equivalent to the category with objects super vector bundles on~$X$ and morphisms are isomorphisms of bundles pulled back to $\pi TX$, i.e., differential forms valued in vector bundle isomorphisms over~$X$. 

With a fixed super vector bundle~$E$, the remaining data of a representation is a section of $\Hom({\sf s}^*E,{\sf t}^*E)$ over ${\rm Mor}(\Path(X\sq G))$. By restricting this section to the submanifold
\beq
\{(0,0)\}\times G\times \pi TX\hookrightarrow \R^{1|1}_{\ge 0}\times G\times \pi TX. \label{eq:Gact}
\eeq
we obtain maps $\rho(g)\colon E\to g^*E$ for each $g\in G$; as we shall see, these determine a twisted equivariant structure. Similarly, we will show restriction to 
\beq
\R^{1|1}_{\ge 0}\times \{e\}\times \pi TX\hookrightarrow \R^{1|1}_{\ge 0}\times G\times \pi TX \label{eq:superconn}
\eeq
determines a super connection compatible with the $G$-action. 
Composition of the morphisms determined by~$\R^{1|1}_{\ge 0}$ and~$G$ commute, so it suffices to understand the representation restricted to~(\ref{eq:Gact}) and~(\ref{eq:superconn}) separately. 

The twisted equivariant structure is easier to see: by the definition of the twist determined by $\beta$, the maps~$\rho(g)$ gotten from the restriction~(\ref{eq:Gact}) are required to compose as
$$
\rho(g)\circ \rho(h)=\beta(g,h)\rho(gh).
$$
Hence,~$\rho$ is precisely a~$\beta$-twisted equivariant structure on~$E$. 


Turning attention to the restriction~(\ref{eq:superconn}), first we observe that the sections of the pullbacks of ${\sf s}^*E$ and ${\sf t}^*E$ to $\R^{1|1}_{\ge 0}\times \pi TX$ can be viewed as sheaves of modules whose underlying vector spaces are both $C^\infty(\R^{1|1}_{\ge 0})\otimes \Omega^\bullet(X;E)$, but where the module structure over functions is different. Specifically, sections of ${\sf s}^*E$ have the obvious module structure over $C^\infty(\R^{1|1}_{\ge 0})\otimes \Omega^\bullet(X)$ and the module structure for ${\sf t}^*E$ is twisted by the algebra map determined by $\Omega^\bullet(X)\to C^\infty(\R^{1|1}_{\ge 0})\otimes \Omega^\bullet(X)$ given by $\omega\mapsto \omega+\theta d\omega$ for $\omega\in \Omega^\bullet(X)$ and $C^\infty(\R^{1|1}_{\ge 0})\cong C^\infty(\R_{\ge 0})[\theta]$ for $\theta$ an odd coordinate. Hence, we identify the restriction~(\ref{eq:superconn}) with a linear endomorphism $R$ of $C^\infty(\R^{1|1}_{\ge 0})\otimes \Omega^\bullet(X;E)$ satisfying
\beq
R(\omega\cdot s)=(\omega+\theta d\omega)\cdot R(s)\label{eq:leib}\quad \omega\in \Omega^\bullet(X), \ s\in C^\infty(\R^{1|1}_{\ge 0})\otimes \Omega^\bullet(X;E). 
\eeq
Compatibility with composition requires that $R$ further satsify
$$
p_1^*R\circ p_2^* R= {\sf c}^*R
$$
for $p_1,p_1,{\sf c}\colon \R^{1|1}_{\ge 0}\times \R_{\ge 0}^{1|1}\times \pi TX\to \R_{\ge 0}^{1|1}\times \pi TX$ the maps associated to the projections and composition. This compatibility determines an operator $\A$ characterized by
$$
(\partial_\theta+\theta\partial_t)R=\A R,
$$
where $(t,\theta)$ are coordinates on $\R^{1|1}_{\ge 0}$. By the existence and uniqueness of solutions to ordinary differential equations on supermanifolds (see~\cite{florin_11}, Section~3.4), we obtain the formula:
$$
R=\exp(-t\A^2+\theta \A),\quad \A\colon \Omega^\bullet(X;E)\to \Omega^\bullet(X;E).
$$
Furthermore, by equation~\ref{eq:leib}, $\A$ satisfies a Leibniz rule,
$$
\A(\omega \cdot s)=d\omega \cdot \A(s)+(-1)^{|\omega|}\omega \cdot \A(s),
$$
and hence $\A$ determines a Quillen super connection on~$E$. 

Compatibility of composition in $\Path(X)$ between morphisms associated to $G$ and morphisms associated to $\R^{1|1}_{\ge 0}$ requires that 
$$
\rho(g)R=R\rho(g),
$$
over $\R^{1|1}_{\ge 0}\times G\times \pi TX$, and so $\A$ defines an \emph{equivariant} super connection on~$E$. Finally, an isomorphism between representations is an isomorphism between vector bundles over~$\pi TX$ (discussed in the first paragraph of the proof) compatible with the operators $\rho(g)$ and $\A$, and hence is an isomorphism of twisted equivariant bundles compatible with the super connection. 
\ep

We make two additional observations about representations of~$\Path(X\sq G)$. First, for a map $X\to Y$ of $G$-manifolds, the induced map $\Path(X\sq G)\to \Path(Y\sq G)$ allows us to pullback representations of $\Path(Y\sq G)$ to representations of $\Path(X\sq G)$, and this pullback agrees with the pullback of twisted equivariant super vector bundles with super connection. Second, there is an evident monoidal structure on the category of representations coming from the direct sum of vector bundles, and we use the notation ${\sf E}_0\oplus{\sf E}_1$ to denote the sum of two twisted representations. 

\subsection{The (Chern) character}

A \emph{super loop} in $X$ is a super path with positive super length (i.e., in $\R_{>0}^{1|1}\subset \R^{1|1}_{\ge 0}$) whose source and target are the same. A \emph{super loop} in $X\sq G$ is a super path with positive super length and a specified isomorphism between its source and target. Hence, the subcategory of~$\Path(X\sq G)$ of energy zero super loops is exactly
$$
\Loop(X\sq G):=\left(\begin{array}{c} 
\R_{> 0}\times \coprod_{g\in G} \pi TX^g\\
{\sf s}\downarrow \downarrow {\sf t}\\
\pi TX
\end{array} \right)\subset \left(\begin{array}{c} 
\R_{\ge 0}^{1|1}\times G\times \SM(\R^{0|1},X)\\
{\sf s}\downarrow \downarrow {\sf t}\\
\SM(\R^{0|1},X)
\end{array} \right)
$$
using the inclusions $\{g\}\times \pi TX^g\subset G\times \pi TX$. 
Geometrically, an $S$-point of $\R_{> 0}\times \coprod_{g\in G} \pi TX^g$ is (1) an $S$-family of super loops $S\times \R^{1|1}/(t\cdot \Z)$ with circumference $t\in \R_{>0}(S)\subset \R^{1|1}_{\ge 0}(S)$; (2) a $G$-bundle over this family with clutching data determined by $g\in G(S)$; and (3) an equivariant map to~$X$ determined by $S\times \R^{0|1}\to X^g$. This map needs to land in $g$-fixed points so that map from the (trivial) $G$-bundle over~$S\times \R^{1|1}$ to $X$ descends to the quotient, i.e., the bundle over $S\times \R^{1|1}/t\Z$ gotten by clutching.

Restriction of a representation of $\Path(X\sq G)$ to $\Loop(X\sq G)$ yields a family of endomorphisms of~${\sf s}^*E={\sf t}^*E$ over $\R_{> 0}\times \coprod_{g\in G} \pi TX^g$. We call the super trace of this endomorphism evaluated at $1\in \R_{>0}$ the \emph{character} of a representation, and we denote it by~$Z({\sf E})$ for a representation~${\sf E}$. Let $Z_g({\sf E})$ denote the restriction of $Z({\sf E})$ to $\pi TX^g$. The data specified by Proposition~\ref{prop:geochar} gives the formula
$$
Z_g({\sf E})={\rm sTr}\left(\exp(-\A^2)\circ \rho(g)\right)\in C^\infty( \pi TX^g)\cong \Omega^\bullet(X),
$$
and so by standard Chern--Weil theory, the character determines an element
$$
Z({\sf E})\in \bigoplus_{g\in G} \Omega^{\ev}_{\rm cl}(X^g). 
$$

The original~$G$-action on~$X$ gives maps~$X^g\to X^{hgh^{-1}}$ for each~$h\in G$, and the character has an important property with respect to this action. First we observe that for a $\beta$-twisted equivariant vector bundle, we have
$$
\rho(hgh^{-1})=\frac{\beta(hgh^{-1},h)}{\beta(h,g)} \rho(h)\circ \rho(g)\circ \rho(h)^{-1},
$$
using the formula for composition and the cocycle condition. We deduce
\beq
{\rm sTr}\left(\exp(-\A^2)\circ \rho(hgh^{-1})\right) =\frac{\beta(hgh^{-1},h)}{\beta(h,g)} {\rm sTr}\left(\exp(-\A^2)\circ \rho(g)\right)\label{eq:partline}
\eeq
from compatibility of $\A$ with the $G$-action and the cyclic property of the super trace. This allows one to view the character~$Z({\sf E})$ as a section of a line bundle over the groupoid $(\coprod_{g\in G} \pi TX^g)\sq G$ that pulls back from the line bundle over $L_0(\pt\sq G)\cong G\sq G$ with cocycle
$$
G\times G\to U(1), (h,g)\mapsto \frac{\beta(hgh^{-1},h)}{\beta(h,g)}\in U(1).
$$
We denote this line bundle by $\mathcal{L}^\beta$ and sections by $\Gamma((\coprod_{g\in G} \pi TX^g)\sq G;\mathcal{L}^\beta)$.

\begin{rmk}
The groupoid $(\coprod_{g\in G} \pi TX^g)\sq G$ is a super version of the \emph{inertia groupoid} of~$X\sq G$, which consists of constant loops in a stack that are permitted to have nontrivial automorphisms connecting their source and target. In our super Euclidean case, the action of loop rotation on $(\coprod_{g\in G} \pi TX^g)\sq G$ is through the de~Rham operator, so we can interpret the closedness of the character as coming from invariance under this rotation action. 
\end{rmk}

\subsection{Chern--Simons forms} 
\begin{defn} \label{def:conc11} A \emph{concordance} between a pair of $\beta$-twisted representations~${\sf E}_0,{\sf E}_1$ of $\Path(X\sq G)$ is a $\beta$-twisted representation $\tilde{\sf E}$ of $\Path((X\sq G)\times \R)$ and isomorphisms $i_0^*\tilde{{\sf E}}\cong {\sf E}_0$, $i_1^*\tilde{{\sf E}}\cong {\sf E}_1$ where $i_0,i_1\colon X\sq G\hookrightarrow X\sq G\times \R$ are the inclusions at $X\sq G\times \{0\}$ and $X\sq G \times \{1\}$. If there is a concordance between ${\sf E}_0$ and ${\sf E}_1$ we call the pair of twisted representations \emph{concordant}. \end{defn} 

Let ${\sf E}_0$ and ${\sf E}_1$ be representations of $\Path(X\sq G)$ corresponding to $\beta$-twisted equivariant super vector bundles with compatible super connections $(E_0,\A^{E_0},\rho_0)$ and $(E_1,\A^{E_1},\rho_1)$. Then ${\sf E}_0$ and ${\sf E}_1$ are concordant if and only if there is an equivariant isomorphism, $\phi\colon E_0\to E_1$. With such an isomorphism fixed, a linear combination of the super connections determines a path $\A(\lambda)$ in the space of super connection with $\A(0)=\A^{E_0}$ and $\A(1)=\phi^*\A^{E_1}$. We can promote this path to a concordance determined by the equivariant vector bundle $p^*E_0$ for $p\colon X\times \R\to X$ with an equivariant structure coming from~$\rho_0$ and super connection $(d\lambda)\partial/\partial\lambda + \A(\lambda)$. The isomorphism $E_1\stackrel{\phi}{\gets} E_0\cong i_1^*(p^*E_0)$ fixes the claimed target of this concordance. 

It will be important to keep track of how the character of a representation changes under a concordance; this is measured by a Chern--Simons form.

\begin{defn} For isomorphic representations $\phi\colon {\sf E}_0\to {\sf E}_1$, define
$$
\CS_g({\sf E}_0,{\sf E}_1,\phi)=\int_{X^g\times I/X^g} Z_g(\tilde{\sf E})\ \ \in \Omega^{\odd}(X^g) \subset C^\infty(\pi TX^g)
$$
where $I=[0,1]$, and the integral is the fiberwise integration of forms and $\tilde{\sf E}$ is the concordance from ${\sf E}_0$ to ${\sf E}_1$ gotten from the linear interpolation between super connections. Let $\CS({\sf E_0},{\sf E}_1,\phi)$ denote the corresponding element of $C^\infty(\coprod_g \pi TX^g)$. 
\end{defn}

The transformation properties of $Z({\sf E})$ under the $G$-action on~$\coprod_g \pi TX^g$ give $\CS({\sf E}_0,{\sf E}_1,\phi)$ the same properties, so we get a section of $\mathcal{L}^\beta$,
$$
\CS({\sf E}_0,{\sf E}_1,\phi)\in \Gamma\left(\Big(\coprod_{g\in G} \pi TX^g\Big)\sq G;\mathcal{L}^\beta\right)\subset \bigoplus_g\Omega^\bullet(X^g)
$$
or more concretely, $\CS({\sf E}_0,{\sf E}_1,\phi)$ is an element of $\Omega^{\odd}(\coprod_g \pi TX^g)$ with specified transformation properties for the $G$-action. 


By the usual arguments, choosing a different path in the space of super connections changes the Chern--Simons form by an exact form. Hence, if we simply know that ${\sf E}_0$ and ${\sf E}_1$ are concordant we have a well-defined odd class
$$
\CS({\sf E}_0,{\sf E}_1)\in \Gamma\left(\Big(\coprod_{g\in G} \pi TX^g\Big)\sq G;\mathcal{L}^\beta\right)/d\Gamma\left(\Big(\coprod_{g\in G} \pi TX^g\Big)\sq G;\mathcal{L}^\beta\right),
$$
i.e., a class in $\Omega^{\rm odd}\Big(\coprod_g X^g\Big)/d\Omega^{\ev}\Big(\coprod_g X^g\Big)$ with properties. 
This class satisfies
$$
Z({\sf E}_1)-Z({\sf E}_0)=d\CS({\sf E}_0,{\sf E}_1),
$$
and so measures how the character of a representation of $\Path(X\sq G)$ changes between concordant field theories. 
 

\subsection{Twisted gauged low-energy effective field theories}

Very roughly, an effective field theory is a finite-dimensional approximation of a (non-effective) field theory. A typical field theory has an infinite-dimensional space of states with an energy filtration coming from the eigenvalues of a Hamiltonian. By imposing an energy cutoff, we obtain a finite-dimensional subspace containing the ``low-energy" states (in particular the kernel of the Hamiltonian). As pertains to the $1|1$-dimensional case, it is precisely this type of finite-dimensional approximation that Atiyah and Singer use to construct the index bundle of a family of Dirac operators. Viewing the square of the Dirac operator as a family of Hamiltonians, our intent is to reinterpret data like an index bundle as an effective field theory. The difference between the character of the finite-dimensional cutoff theory and the infinite-dimensional one defines a sort of Chern--Simons form that we take as data.

\begin{defn} A \emph{$1|1$-dimensional, $\beta$-twisted, gauged, low-energy effective field theory over~$X$} is a $\beta$-twisted representation ${\sf E}$ of $\Path(X\sq G)$ and an odd element
$$
\eta \in \Gamma\left(\Big(\coprod_{g\in G} \pi TX^g\Big)\sq G;\mathcal{L}^\beta\right)/d\Gamma\left(\Big(\coprod_{g\in G} \pi TX^g\Big)\sq G;\mathcal{L}^\beta\right),
$$ 
i.e., an equivalence class~$\eta\in \Omega^{\rm odd}\Big(\coprod_g X^g\Big)/d\Omega^{\ev}\Big(\coprod_g X^g\Big)$ with specified transformation properties under the $G$-action. Low-energy effective field theories $({\sf E}_0,\eta_0)$ and $({\sf E}_1,\eta_1)$ are \emph{isomorphic} if their representations are concordant and there is an equality $\eta_0=\eta_1+\CS({\sf E}_0,{\sf E}_1)$. Define a monoidal structure $\oplus$ on this category by 
$$
({\sf E}_0,\eta_0)\oplus ({\sf E}_1,\eta_1)=({\sf E}_0\oplus {\sf E}_1,\eta_0+\eta_1),
$$
and denote this monoidal category by $1|1\EFT^\beta(X\sq G)$. 
\end{defn}

We will often drop the adjectives ``$1|1$-dimensional, $\beta$-twisted, gauged, low-energy" and ``over $X$" when they are implied by context. 

\begin{rmk} In a sense, these effective field theories are of low-energy in two different ways: (1) they only consider energy zero paths in~$X$ and (2) they are a low-energy approximation of a possibly infinite-dimensional theory. 
\end{rmk}

\begin{rmk} The \emph{partiton function} of an effective field theory is the differential form
$$
Z({\sf E},\eta)=Z({\sf E})+d\eta\in \Gamma\left(\Big(\coprod_{g\in G} \pi TX^g\Big)\sq G;\mathcal{L}^\beta\right)\subset \bigoplus_{g\in G}\Omega^{\ev}_{\rm cl}(X^g). 
$$
\end{rmk}

It remains to explain a notion of \emph{stable} isomorphism; physically this mediates between effective field theories that arise from different choices of cutoff energy. 

%

\begin{defn}
A twisted representation $\Path(X\sq G)\to \Vect$ is \emph{stably trivial} if it is isomorphic to one whose super vector bundle is of the form $V\oplus \pi V$ for $V$ an ordinary (purely even) twisted equivariant vector bundle on~$X$  whose super connection comes from an ordinary equivariant connection $\nabla$ on~$V$. We denote the trivial representation determined by $(V,\nabla)$ by $\epsilon_V$. 
\end{defn}

\begin{defn} Two low-energy effective field theories $({\sf E}_0,\eta_0)$, $({\sf E}_1,\eta_1)$ are \emph{stably isomorphic} if $\eta_0=\eta_1$ and there is an isomorphism of representations ${\sf E}_0\cong E_1\oplus \epsilon_V.$
\end{defn}

\begin{rmk} Analogous to the situation of an index bundle, the difference between two cutoff theories is a subspace on which~$\A_0$ is invertible; if~$\A_0$ is self-adjoint for some metric, $\A_0^2$~has strictly positive eigenvalues corresponding to positive energies. But a representation with $\A_0$ invertible is concordant to a trivial representation since the even and odd parts of their super vector bundles are isomorphic via~$\A_0$. Hence, in the category of effective field theories we get an isomorphism to a trivial representation with a class~$\eta$ coming from the concordance to the trivial representation. \end{rmk}

\begin{rmk} One can also view trivial field theories as defining ones that limit to the zero theory under the \emph{renormalization group} (RG) flow. To spell this out, define a super connection on $V\oplus \pi V$ whose degree zero piece $\A_0$ comes from the identity map $V\to V$, viewed as an odd endomorphism of $V\oplus \pi V$, and whose degree 1 piece $\A_1$ comes from an ordinary connection on~$V$. The RG flow on field theories is induced from dilating the super length of super paths; in our case this leads to the $\R_{>0}$-family of super connections $\A(\lambda)=\lambda\A_0+\A_1$ for $\lambda\in \R_{>0}$. In the limit~$\lambda\to\infty$, $\exp(-t\A(\lambda)^2+\theta\A(\lambda))$ is the zero operator on~$V\oplus \pi V$ which is the sense in which this theory limits to zero. Furthermore, the Chern--Simons form for a concordance from~$\A(0)$ to~$\A(\lambda)$ is zero for all~$\lambda$, so there is no affect on the character. 
\end{rmk}

\begin{rmk} Stolz and Teichner have defined a category of $1|1$-Euclidean field theories whose objects are functors from a bordism category of $1|1$-Euclidean manifolds over~$X$ to topological vector spaces. Results of Dumitrescu~\cite{florin_11} on super parallel transport along a super connection allow us to extend a representation of $\Path(X)$ to a field theory in the sense of Stolz and Teichner~\cite{ST11}. Conversely, for Stolz--Teichner field theories valued in finite-dimensional vector spaces, restriction to energy zero super paths defines a representation in our sense. In this way, one can view representations of $\Path(X)$ as encoding the subcategory of Stolz and Teichner's $1|1$-Euclidean field theories that take values in finite-dimensional vector spaces. 
\label{rmk:ST}\end{rmk}


\section{Low-energy effective field theories and differential K-theory}\label{sec:diffK}

\subsection{Twisted equivariant differential K-theory}

We follow the model for equivariant differential K-theory developed by Ortiz~\cite{Ortiz} Section~3.3, which in turn follows the general description of a differential cohomology theory put forward by Hopkins and Singer~\cite{hopsing}. The starting point is a description of $K_G^\beta(X)\otimes \C$. 

\begin{thm}[\cite{AdemRuan} Theorem~7.4] Let $\beta\colon G\times G\to U(1)$ be a normalized 2-cocycle. For a manifold $X$ with $G$-action, the $\beta$-twisted $G$-equivariant ${\rm K}$-theory of~$X$ with complex coefficients can be computed as
$$
K_G^{\beta}(X)\otimes \C \cong\bigoplus_{[g]} ({\rm H}_{\rm dR}^{\ev}(X^g)\otimes \chi^\beta_g)^{C_G(g)}
$$
where the sum ranges over conjugacy classes of $g\in G$, $C_G(g)$ denotes the centralizer of $g$, and $\chi^\beta_g$ denotes the 1-dimensional representation of $C_G(g)$ given by $h\mapsto \beta(h,g)\beta(g,h)^{-1}$. \label{thm:KtensorC}
\end{thm}

\begin{defn} \label{def:DiffK}
Let $\widehat{\K}_G^\beta(X)$ be the free abelian group on generators $(V,\nabla,\eta)$ for $V$ a $\Z/2$-graded $\beta$-twisted equivariant vector bundle on $X$ with (grading-preserving) connection $\nabla$ and $\eta\in \bigoplus_{[g]} (\Omega^{\odd}(X^g)\otimes \chi^\beta_g/d\Omega^{\ev}(X^g)\otimes \chi^\beta_g)^{C_G(g)}$ is an equivalence class of odd forms subject to the relations
\begin{enumerate}
\item If $V$ and $V'$ are isomorphic, then $(V,\nabla,\eta)\sim (V',\nabla',\eta+ \CS(\nabla',\nabla))$;
\item $(V,\nabla,\eta)+(V',\nabla',\eta')\sim (V\oplus V',\nabla\oplus \nabla',\eta+\eta')$; and
\item $(V,\nabla,\eta)+(\pi V,\nabla,-\eta)\sim 0$.
\end{enumerate}
\end{defn}

\subsection{Proof of the main theorem} The cocycle for the line bundle $\mathcal{L}^\beta$ under the equivalence of stacks induced by the inclusion
$$
\coprod_{[g]} \Big(\pi TX^g\sq C_G(g)\Big)\subset \Big(\coprod_g \pi TX^g\Big)\sq G
$$
is $(h,g)\mapsto \beta(h,g)\beta(g,h)^{-1}$. From this we have an evident inclusion of cocycles in Definition~\ref{def:DiffK} into the category of effective field theories: $(V,\nabla)$ determine a $\beta$-twisted representation of $\Path(X\sq G)$, and~$\eta$ encodes a section of $\mathcal{L}^\beta$ under the equivalence above. From relation (2) in Definition~\ref{def:DiffK}, the abelian group structure on differential cocycles is compatible with the monoid structure on isomorphism classes of effective field theories. 

Furthermore, every effective field theory is isomorphic to one for which the underlying twisted representation comes from an grading-preserving connection: the space of super connections is affine, so we can always choose a concordance of the representation and modify the differential form data accordingly. Thus, we have given a bijection between isomorphism classes of effective field theories and the free abelian group on triples $(V,\nabla,\eta)$ modulo the relations~(1) and~(2) in Definition~\ref{def:DiffK}.

Under this bijection, the relation of stable equivalence is precisely the relation~(3) in Definition~\ref{def:DiffK}, and hence we have a bijection,
$$
\widehat{\K}^\beta_G(X)\cong 1|1\EFT^\beta(X\sq G)/{\sim}
$$
and the theorem is proved. We observe that under this isomorphism, the curvature map in differential K-theory agrees with the partition function of the associated effective field theory. 

%
%
%
%

\bibliographystyle{amsalpha}
\bibliography{references}

\end{document}